\documentclass[12pt]{article}
\usepackage{amsbsy,amsfonts,amsmath,amssymb,amsthm,euscript,enumerate,ascmac}
\usepackage{graphicx,color}
\usepackage{url}
\usepackage{hyperref}
\hypersetup{colorlinks,citecolor=blue,linkcolor=blue,urlcolor=black}
\usepackage{bbm}
\newcommand{\indic}{\mathbbm{1}}

\usepackage{subcaption}
\definecolor{mygreen}{RGB}{0,128,0}

\allowdisplaybreaks[1]

\theoremstyle{plain}
\numberwithin{equation}{section}
\newtheorem{thm}{Theorem}[section]

\newtheorem{dfn}[thm]{Definition}

\newtheorem{prp}[thm]{Proposition}
\newtheorem{rmk}[thm]{Remark}

\newcommand{\dpst}{\displaystyle}
\newcommand{\EG}{E^{\sss\rm G}}
\newcommand{\ESCA}{E^{\sss\rm SCA}}

\newcommand{\GS}{{\sss\rm GS}}

\newcommand{\ind}[1]{\indic_{\{#1\}}}
\newcommand{\lbeq}[1]{\label{eq:#1}}
\newcommand{\mE}{{\mathbb E}}

\newcommand{\nn}{\nonumber}

\newcommand{\PG}{P^{\sss\rm G}}
\newcommand{\PSCA}{P^{\sss\rm SCA}}
\newcommand{\piG}{\pi^{\sss\rm G}}
\newcommand{\piSCA}{\pi^{\sss\rm SCA}}
\newcommand{\Proof}[1]{\paragraph{\it #1}}
\newcommand{\QED}{\hspace*{\fill}\rule{7pt}{7pt}\smallskip}
\newcommand{\refeq}[1]{(\ref{eq:#1})}
\newcommand{\sss}{\scriptscriptstyle}

\newcommand{\tvdist}[1]{\|#1\|_{\sss\text{TV}}}

\newcommand{\vep}{\varepsilon}

\newcommand{\wG}{w^{\sss\rm G}}
\newcommand{\wSCA}{w^{\sss\rm SCA}}

\newcommand{\sigmavec}  {\boldsymbol{\sigma}}
\newcommand{\tauvec}  {\boldsymbol{\tau}}

\title{Finding optimal solutions by stochastic cellular automata}
\author{
Satoshi Handa\footnote{Graduate School of Mathematics, 
Hokkaido University. (Current affiliation: Fujitsu Laboratories Ltd.)}\\
Katsuhiro Kamakura\footnote{Graduate School of Mathematics, Hokkaido University.}\\
Yoshinori Kamijima\footnotemark[2]\\
Akira Sakai\footnote{Faculty of Science, Hokkaido University. 
\url{https://orcid.org/0000-0003-0943-7842}}
}

\begin{document}
\maketitle

\begin{abstract}
Finding a ground state of a given Hamiltonian is an important but hard problem. 
One of the potential methods is to use a Markov chain Monte Carlo (MCMC) to 
sample the Gibbs distribution whose highest peaks correspond to the ground 
states.  In this short paper, we use stochastic cellular  automata (SCA) and 
see if it is possible to find a ground state faster than the conventional 
MCMCs, such as the Glauber dynamics.  We show that, if the temperature is 
sufficiently high, it is possible for SCA to have more spin-flips per update in 
average than Glauber and, at the same time, to have an equilibrium distribution 
``close" to the one for Glauber, i.e., the Gibbs distribution.  During the 
course, we also propose a new way to characterize how close a probability 
measure is to the target Gibbs.
\end{abstract}

\section{Introduction}
There are many occasions in real life when we have to quickly choose one 
among extremely many options.  In addition, we want our choice to be optimal in 
a certain sense.  Such optimization problems are ubiquitous and possibly quite 
hard to solve them fast. 
In particular, NP-hard problems cannot be solved in polynomial time~\cite{gj79}.

One approach to find an optimal solution to a given problem is to translate it 
into an Ising model (see, e.g., \cite{l14} for a list of examples of such 
mappings) and find its ground state that corresponds to an optimal solution.  
A standard method to do so is to use a Markov chain Monte Carlo (MCMC) to 
sample the Gibbs distribution
$\piG_\beta(\sigmavec)\propto e^{-\beta H(\sigmavec)}$, where $\beta\ge0$ is 
the inverse temperature and $H(\sigmavec)$ is the corresponding Ising 
Hamiltonian.  The ground states are the spin configurations at which the 
Gibbs distribution takes on its highest peaks.

There are several MCMCs that generate the Gibbs distribution.  In most of the 
popular ones, such as the Glauber dynamics~\cite{g63}, the number of spin-flips 
per update is at most one, which makes them inefficient.  We have longed for 
faster MCMCs.

We consider probabilistic cellular automata, or PCA for short 
\cite{dss12,st18}.  Since the term PCA has already been used for long as an 
abbreviation for principal component analysis in statistics, we should rather 
call it stochastic cellular automata (SCA).  It is an MCMC based on independent 
multi-spin flip dynamics, and therefore it is potentially faster than the 
standard single-spin flip MCMCs.  

However, since the SCA equilibrium $\piSCA_{\beta,q}$ is different from 
$\piG_\beta$, we cannot naively use it as a Gibbs sampler to search for the 
ground states.  Only when the pinning parameter $q\ge0$ goes to infinity, 
the total variation distance $\tvdist{\piSCA_{\beta,q}-\piG_\beta}$ 
goes to zero.  The downside of this limit is the slowdown effect on SCA: each 
spin is more likely to stick to the current state for larger values of $q$ 
(this is why we call it the pinning parameter).

In this paper, we investigate SCA and try to find by simple arithmetic an 
interval of $q$ that is compatible with good approximation to Gibbs and, at the 
same time, faster than Glauber.  
In Section~\ref{s:definition}, we define relevant observables, such as the 
Hamiltonian, the transition probabilities for Glauber and SCA, and their 
equilibrium distributions.
In Section~\ref{s:morespins}, we show an upper bound on $q$ below which SCA is 
faster than Glauber in terms of the expected number of spin-flips per update, 
not in terms of the mixing time.  In Section~\ref{s:close}, we show a lower 
bound on $q$ above which the SCA equilibrium $\piSCA_{\beta,q}$ is close 
to the target Gibbs $\piG_\beta$ in the sense of order-preservation (see 
\refeq{closedef} for the precise definition), rather than in the total variation 
distance.  Thanks to this new notion of closeness, it becomes tractable to 
obtain a quantitative estimate on the difference between the two equilibrium 
measures.  This is crucial especially when we apply the theory to real-life 
situations or numerical simulations.
In the last section, Section~\ref{s:concluding}, we give a conclusion obtained 
from the results in Sections~\ref{s:morespins}--\ref{s:close}, which can be 
stated in short words as follows.

\begin{thm}
If $\beta$ is sufficiently small, then we can define SCA for which the following 
both hold at the same time:
\begin{itemize}
\item
SCA is faster than Glauber in terms of the expected number of spin-flips per 
update.
\item
If the SCA equilibrium $\piSCA_{\beta,q}$ takes on its highest peak at 
$\sigmavec$, then $\sigmavec$ is close to the ground states of $H$ 
up to a certain error (see Theorem~\ref{thm:summary} for the precise statement).
\end{itemize}
\end{thm}

At the end of this paper, we also discuss implications of the above 
theorem and an idea to find the ground states without using any cooling 
methods.

\section{Definition}\label{s:definition}
Given a finite graph $G=(V,E)$ with no multi- or self-edges, spin-spin 
couplings $\{J_{x,y}\}_{\{x,y\}\in E}$ (with $J_{x,y}=0$ if $\{x,y\}\notin E$) 
and local magnetic fields $\{h_x\}_{x\in V}$, we define the Hamiltonian for 
the spin configuration 
$\sigmavec=\{\sigma_x\}_{x\in V}\in\{\pm1\}^V$ as
\begin{align}
H(\sigmavec)=-\sum_{\{x,y\}\in E}J_{x,y}\sigma_x\sigma_y-\sum_{x\in V}h_x
 \sigma_x=-\frac12\sum_{x,y\in V}J_{x,y}\sigma_x\sigma_y-\sum_{x\in V}h_x
 \sigma_x.
\end{align}
Let $\wG_\beta$ and $\piG_\beta$ be the Boltzmann weight and the Gibbs 
distribution, respectively, at the inverse temperature $\beta\ge0$:
\begin{align}
\wG_\beta(\sigmavec)=e^{-\beta H(\sigmavec)},&&
\piG_\beta(\sigmavec)=\frac{\wG_\beta(\sigmavec)}{\sum_{\sigmavec}
 \wG_\beta(\sigmavec)}.
\end{align}
For $\sigmavec\in\{\pm1\}$, $I\subset V$ and $x\in V$, we let
\begin{align}\lbeq{sigmavec^I}
(\sigmavec^I)_y=
 \begin{cases}
 \sigma_y&[y\notin I],\\
 -\sigma_y\quad&[y\in I],
 \end{cases}&&
\sigmavec^x=\sigmavec^{\{x\}}.
\end{align}
It is known that $\piG_\beta$ is the equilibrium distribution for the Glauber 
dynamics whose transition probability is defined as
\begin{align}\lbeq{PGdef}
\PG_\beta(\sigmavec,\tauvec)=
 \begin{cases}
 \dpst\frac1{|V|}\frac{\wG_\beta(\sigmavec^x)}{\wG_\beta(\sigmavec)+\wG_\beta
  (\sigmavec^x)}\quad&[\tauvec=\sigmavec^x],\\[1pc]
 \dpst1-\sum_{x\in V}\PG_\beta(\sigmavec,\sigmavec^x)&[\tauvec=\sigmavec],\\
 0&[\text{otherwise}].
 \end{cases}
\end{align}
By introducing the cavity fields $\tilde h_x(\sigmavec)$, defined as
\begin{align}\lbeq{cavity}
\tilde h_x(\sigmavec)=\sum_{y\in V}J_{x,y}\sigma_y+h_x,
\end{align}
the Glauber transition probability $\PG_\beta(\sigmavec,\sigmavec^x)$ can be 
written as
\begin{align}\lbeq{PGrewr}
\PG_\beta(\sigmavec,\sigmavec^x)=\frac1{|V|}\frac{e^{-\beta\tilde h_x
 (\sigmavec)\sigma_x}}{2\cosh(\beta\tilde h_x(\sigmavec))}.
\end{align}

Next we define SCA.  First we let
\begin{align}
\tilde H(\sigmavec,\tauvec)=-\frac12\sum_{x,y\in V}J_{x,y}\sigma_x\tau_y-\frac12\sum_{x
 \in V}h_x(\sigma_x+\tau_x).
\end{align}
Notice that
\begin{align}\lbeq{symmetry}
\tilde H(\sigmavec,\tauvec)=\tilde H(\tauvec,\sigmavec),&&
\tilde H(\sigmavec,\sigmavec)=H(\sigmavec).
\end{align}
Let
\begin{align}
\wSCA_{\beta,q}(\sigmavec)&=\sum_{\tauvec}\exp\bigg(-\beta\tilde
 H(\sigmavec,\tauvec)+q\sum_{x\in V}\sigma_x\tau_x\bigg),\\
\piSCA_{\beta,q}(\sigmavec)&=\frac{\wSCA_{\beta,q}(\sigmavec)}{\sum_{\sigmavec}
 \wSCA_{\beta,q}(\sigmavec)}.
\end{align}
It is easy to see from \refeq{symmetry} 
that $\piSCA_{\beta,q}$ satisfies the detailed-balance condition with the 
following transition probability:
\begin{align}\lbeq{PSCAdef}
\PSCA_{\beta,q}(\sigmavec,\tauvec)=\frac1{\wSCA_{\beta,q}(\sigmavec)}\exp
 \bigg(-\beta\tilde H(\sigmavec,\tauvec)+q\sum_{x\in V}\sigma_x\tau_x\bigg).
\end{align}
By using the cavity fields \refeq{cavity}, we can rewrite
$\PSCA_{\beta,q}(\sigmavec,\tauvec)$ as
\begin{align}\lbeq{PSCAredef}
\PSCA_{\beta,q}(\sigmavec,\tauvec)=\prod_{x\in V}\frac{e^{(\frac12\beta\tilde
 h_x(\sigmavec)+q\sigma_x)\tau_x}}{2\cosh(\frac12\beta\tilde h_x(\sigmavec)+q
 \sigma_x)},
\end{align}
which implies independent updating of all spins, and that SCA can jump from one 
spin configuration to any other in a single step.

We note that, since the underlying graph $G$ is finite, the total variation 
distance $\tvdist{\piSCA_{\beta,q}-\piG_\beta}\equiv\frac12\sum_{\sigmavec}
|\piSCA_{\beta,q}(\sigmavec)-\piG_\beta(\sigmavec)|$ gets smaller as 
$q\uparrow\infty$, while the transition probability 
$\PSCA_{\beta,q}(\sigmavec,\tauvec)$ for $\tauvec\ne\sigmavec$ also gets 
smaller as $q\uparrow\infty$ due to the pinning term 
$q\sum_{x\in V}\sigma_x\tau_x$ in \refeq{PSCAdef}, potentially resulting in 
slower convergence.  Our goal is to find an interval of $q$ that guarantees 
faster convergence to $\piSCA_{\beta,q}$ which is close enough to Gibbs in a 
certain sense.

\section{A sufficient condition for more spin-flips}\label{s:morespins}
Given two spin configurations $\sigmavec,\tauvec\in\{\pm1\}^V$, we let 
$D_{\sigmavec,\tauvec}$ be the set of vertices at which $\sigmavec$ and 
$\tauvec$ take on different values:
\begin{align}
D_{\sigmavec,\tauvec}=\{x\in V:\sigma_x\ne\tau_x\}.
\end{align}
We want to compare the expected number 
$E^*[|D_{\sigmavec,X_{\sigmavec}}|]\equiv\sum_{\tauvec}|D_{\sigmavec,
\tauvec}|P^*(\sigmavec,\tauvec)$ of spin-flips per update, 
where $X_{\sigmavec}$ is a $\{\pm1\}^V$-valued random variable whose law 
is $P^*(\sigmavec,\cdot)$,
between Glauber ($P^*=\PG_\beta$) and SCA ($P^*=\PSCA_{\beta,q}$), 
uniformly in the initial configuration $\sigmavec$.  Here is a sufficient 
condition for the latter to be larger than the former.
\begin{prp}\label{prp:morespins}
Let 
\begin{align}\lbeq{Kdef}
\bar K=\max\big\{|\tilde h_x(\sigmavec)|:x\in V,~\sigmavec\in\{\pm1\}^V\big\}.
\end{align}
Then, 
$\EG_\beta[|D_{\sigmavec,X_{\sigmavec}}|]\le\ESCA_{\beta,q}[|D_{\sigmavec,X_{\sigmavec}}|]$ 
holds for any $\sigmavec$ if 
\begin{align}\lbeq{morespins}
2q\le\log|V|-\beta\bar K,
\end{align}
which is positive if $\beta$ is sufficiently small, depending on $G$, 
$\{J_{x,y}\}_{\{x, y\}\in E}$ and $\{h_x\}_{x\in V}$.
\end{prp}

\begin{rmk}\label{rem:morespins}
{\rm
\begin{enumerate}[(i)]
\item
We note that the inequality 
$\EG_\beta[|D_{\sigmavec,X_{\sigmavec}}|]\le\ESCA_{\beta,q}[|D_{\sigmavec,
X_{\sigmavec}}|]$ does not necessarily mean that the SCA convergence to 
equilibrium is faster than that of Glauber.  To prove faster convergence in 
SCA, we must compare their spectral gaps or the mixing times \cite{lpw09}.  
This is under investigation with Bruno Kimura.
\item
In \cite{dss12}, the expected number of spin-flips per update is claimed 
to be of order $|V|e^{-2q}$.  However, this is a bit misleading, as explained 
now.  First, we recall \refeq{PSCAredef}.  Notice that
\begin{align}
\frac{e^{(\frac\beta2\tilde h_x(\sigmavec)+q\sigma_x)\tau_x}}{2\cosh(\frac\beta2
 \tilde h_x(\sigmavec)+q\sigma_x)}
&=\ind{x\in D_{\sigmavec,\tauvec}}\frac{e^{-\frac\beta2\tilde h_x(\sigmavec)
 \sigma_x-q}}{2\cosh(\frac\beta2\tilde h_x(\sigmavec)+q\sigma_x)}\nn\\
&\quad+\ind{x\in V\setminus D_{\sigmavec,\tauvec}}\frac{e^{\frac\beta2\tilde h_x
 (\sigmavec)\sigma_x+q}}{2\cosh(\frac\beta2\tilde h_x(\sigmavec)+q\sigma_x)}.
\end{align}
Isolating the $q$-dependence, we can rewrite the first factor on the right as 
\begin{align}
\frac{e^{-\frac\beta2\tilde h_x(\sigmavec)\sigma_x-q}}{2\cosh(\frac\beta2\tilde
 h_x(\sigmavec)+q\sigma_x)}=\underbrace{\frac{e^{-q}\cosh(\frac\beta2\tilde
 h_x(\sigmavec))}{\cosh(\frac\beta2\tilde h_x(\sigmavec)+q\sigma_x)}}_{\equiv\,
 \vep_x(\sigmavec)}\,\underbrace{\frac{e^{-\frac\beta2\tilde h_x(\sigmavec)
 \sigma_x}}{2\cosh(\frac\beta2\tilde h_x(\sigmavec))}}_{\equiv\,p_x(\sigmavec)},
\end{align}
and the second factor on the right as $1-\vep_x(\sigmavec)p_x(\sigmavec)$.  
As a result, we obtain
\begin{align}
\PSCA_{\beta,q}(\sigmavec,\tauvec)=\prod_{x\in D_{\sigmavec,\tauvec}}\big(\vep_x
 (\sigmavec)p_x(\sigmavec)\big)\prod_{y\in V\setminus D_{\sigmavec,\tauvec}}
 \big(1-\vep_y(\sigmavec)p_y(\sigmavec)\big).
\end{align}
Suppose that $\vep_x(\sigmavec)$ is independent of $x$ and $\sigmavec$, which 
is of course untrue, and simply denote it by $\vep=O(e^{-2q})$.  Then we can 
rewrite $\PSCA_{\beta,q}(\sigmavec,\tauvec)$ as
\begin{align}
\PSCA_{\beta,q}(\sigmavec,\tauvec)
&\simeq\prod_{x\in D_{\sigmavec,\tauvec}}\big(
 \vep p_x(\sigmavec)\big)\prod_{y\in V\setminus D_{\sigmavec,\tauvec}}\Big((1
 -\vep)+\vep\big(1-p_y(\sigmavec)\big)\Big)\nn\\
&=\prod_{x\in D_{\sigmavec,\tauvec}}\big(\vep p_x(\sigmavec)\big)\sum_{S:
 D_{\sigmavec,\tauvec}\subset S\subset V}(1-\vep)^{|V\setminus S|}\prod_{y\in
 S\setminus D_{\sigmavec,\tauvec}}\Big(\vep\big(1-p_y(\sigmavec)\big)\Big)\nn\\
&=\sum_{S:D_{\sigmavec,\tauvec}\subset S\subset V}\vep^{|S|}(1-\vep)^{|V
 \setminus S|}\prod_{x\in D_{\sigmavec,\tauvec}}p_x(\sigmavec)\prod_{y\in S
 \setminus D_{\sigmavec,\tauvec}}\big(1-p_y(\sigmavec)\big).
\end{align}
This implies that the transition from $\sigmavec$ to $\tauvec$ can be seen as 
determining the binomial subset $D_{\sigmavec,\tauvec}\subset S\subset V$ with 
parameter $\vep$ and then changing each spin at $x\in D_{\sigmavec,\tauvec}$ 
with probability $p_x(\sigmavec)$.  Therefore, $|V|\vep$ could be much larger 
than the actual expected number of spin-flips per update.
\end{enumerate}
}
\end{rmk}

\Proof{Proof of Proposition~\ref{prp:morespins}.}
For Glauber, since there is at most one spin-flip, we have
\begin{align}
\EG_\beta[|D_{\sigmavec,X_{\sigmavec}}|]=\sum_{x\in V}\PG_\beta(\sigmavec,
 \sigmavec^x)
&=\frac1{|V|}\sum_{x\in V}\frac{e^{-\beta\tilde h_x(\sigmavec)\sigma_x}}{2
 \cosh(\beta\tilde h_x(\sigmavec))}\nn\\
&=\frac1{|V|}\sum_{x\in V}\big(e^{2\beta\tilde h_x(\sigmavec)\sigma_x}+1
 \big)^{-1}.
\end{align}
On the other hand, since 
$|D_{\sigmavec,\tauvec}|=\sum_{x\in V}\ind{\sigma_x\ne\tau_x}$, we have
\begin{align}
\ESCA_{\beta,q}[|D_{\sigmavec,X_{\sigmavec}}|]=\sum_{x\in V}\sum_{\tauvec:\tau_x
 \ne\sigma_x}\PSCA_{\beta,q}(\sigmavec,\tauvec)&=\sum_{x\in V}\frac{e^{-\frac12
 \beta\tilde h_x(\sigmavec)\sigma_x-q}}{2\cosh(\frac12\beta\tilde h_x(\sigmavec)
 \sigma_x+q)}\nn\\
&=\sum_{x\in V}\big(e^{\beta\tilde h_x(\sigmavec)\sigma_x+2q}+1\big)^{-1}.
\end{align}
A sufficient condition for $\EG_\beta[|D_{\sigmavec,X_{\sigmavec}}|]\le
\ESCA_{\beta,q}[|D_{\sigmavec,X_{\sigmavec}}|]$ is the term-wise inequality
\begin{align}
|V|\big(e^{2\beta\tilde h_x(\sigmavec)\sigma_x}+1\big)\ge e^{\beta\tilde
 h_x(\sigmavec)\sigma_x+2q}+1,
\end{align}
which immediately follows from \refeq{morespins} (by ignoring 1 in both sides of 
the above inequality).
\QED

\section{A sufficient condition for being close to Gibbs}\label{s:close}
To characterize how close a probability measure $\mu$ on $\{\pm1\}^V$ is to 
the target Gibbs $\piG_\beta$, the total variation distance 
$\tvdist{\mu-\piG_\beta}$ has been a standard norm in the 
literature.  If we try to find the ground states by using this norm with 
$\mu=\piSCA_{\beta,q}$, then we have to take $q$ so large that 
$\piSCA_{\beta,q}$ is to some extent uniformly close to $\piG_\beta$.  Such a 
large pinning parameter $q$, however, may not be compatible with the inequality 
\refeq{morespins}.  This is the downside of using total variation.

\begin{figure}[t]
\begin{align*}
\includegraphics[scale=.35]{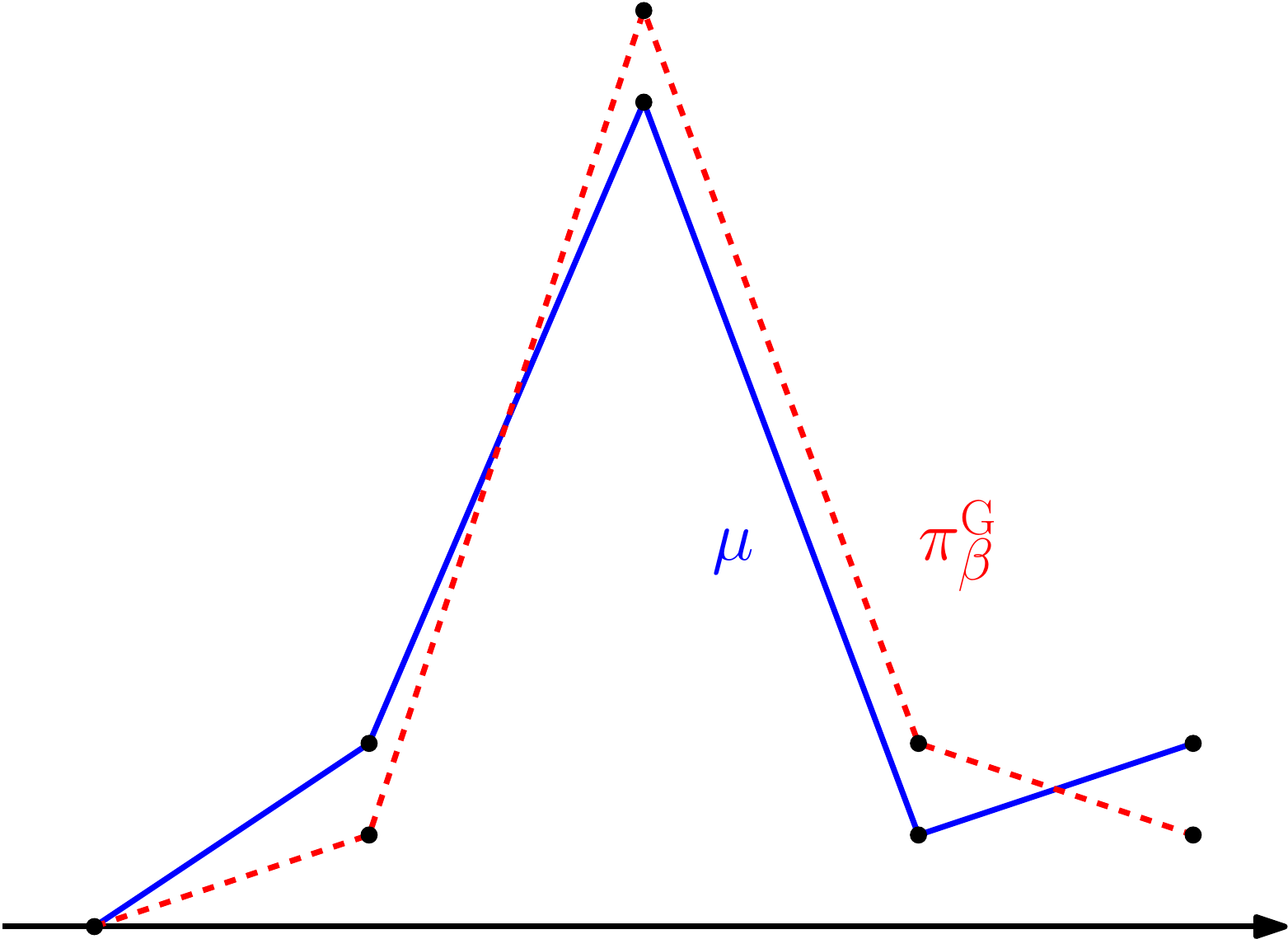}&&&&
\includegraphics[scale=.35]{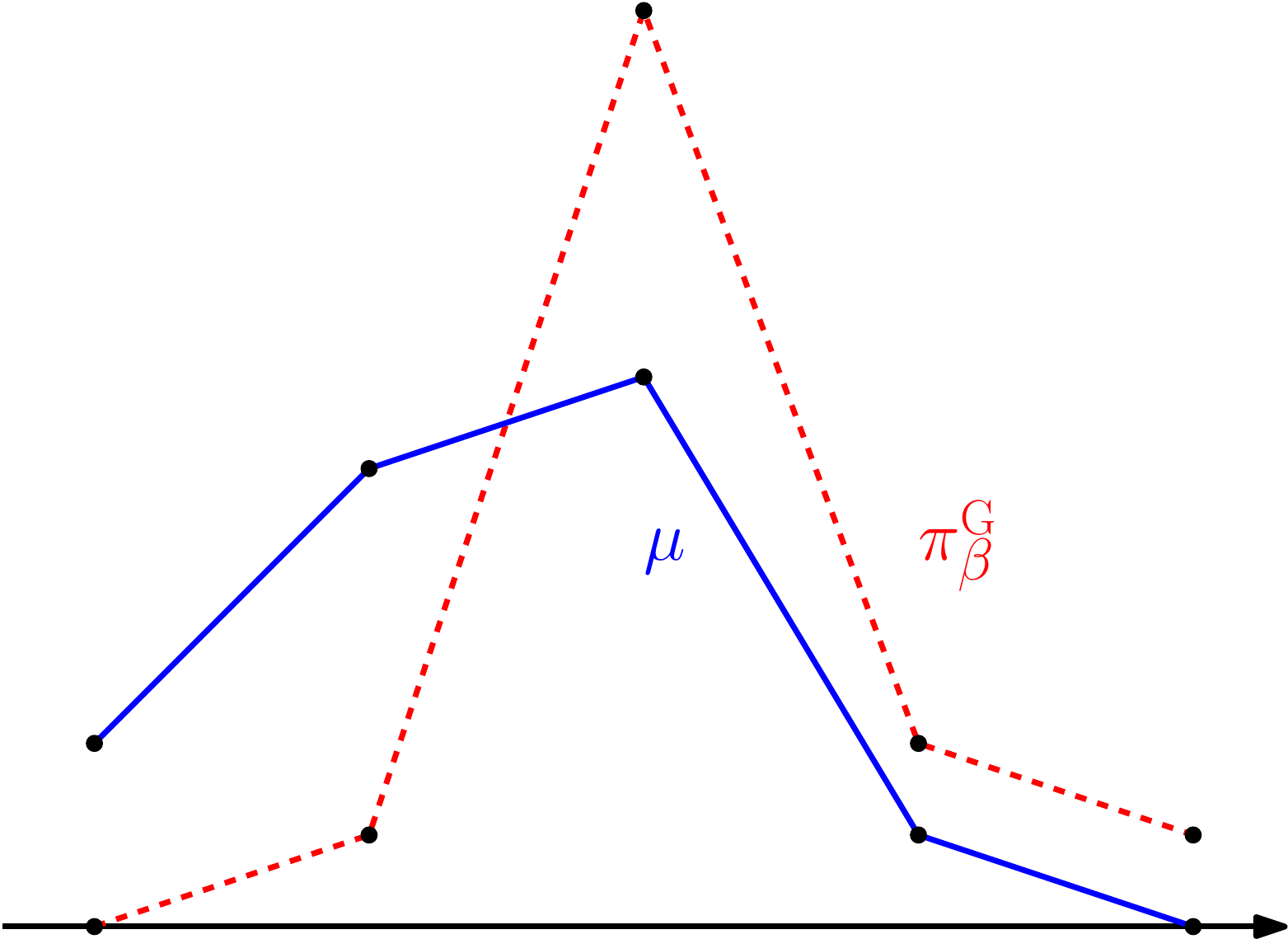}
\end{align*}
\caption{\label{fig:dist}
Comparison of two distributions $\mu$ and $\piG_\beta$ in two cases: 
either the total variation distance $\tvdist{\mu-\piG_\beta}$ is 
small (left) or is not (right).  We can correctly find the ground state of 
$\piG_\beta$ by $\mu$ in both cases.}
\end{figure}

We realize that, to find the ground states $\{\sigmavec_\GS\}$ by another 
distribution $\mu$, it does not have to be close to $\piG_\beta$ in total 
variation, but we rather want it to take the highest peaks among 
$\{\sigmavec_\GS\}$ (see Figure~\ref{fig:dist}).

Inspired by this observation, we introduce the following new notion of 
closeness to the target Gibbs, which is also used in the stability analysis in 
\cite{nstu??}.

\begin{dfn}\label{def:close}
{\rm
Let $\vep>0$ and define the range of $H$ as
\begin{align}
R_H=\max_{\sigmavec,\tauvec\in\{\pm1\}^V}|H(\sigmavec)-H(\tauvec)|.
\end{align}
We say that a probability measure $\mu$ on $\{\pm1\}^V$ is \emph{$\vep$-close 
to the target Gibbs $\piG_\beta$ in the sense of order-preservation} if 
\begin{align}\lbeq{closedef}
\mu(\sigmavec)\ge\mu(\tauvec)
\quad\Rightarrow\quad
H(\sigmavec)\le H(\tauvec)+\vep R_H.
\end{align}
}
\end{dfn}

\bigskip

\begin{figure}[t]
  \begin{subfigure}{0.48\linewidth}
    \includegraphics[width=\linewidth]{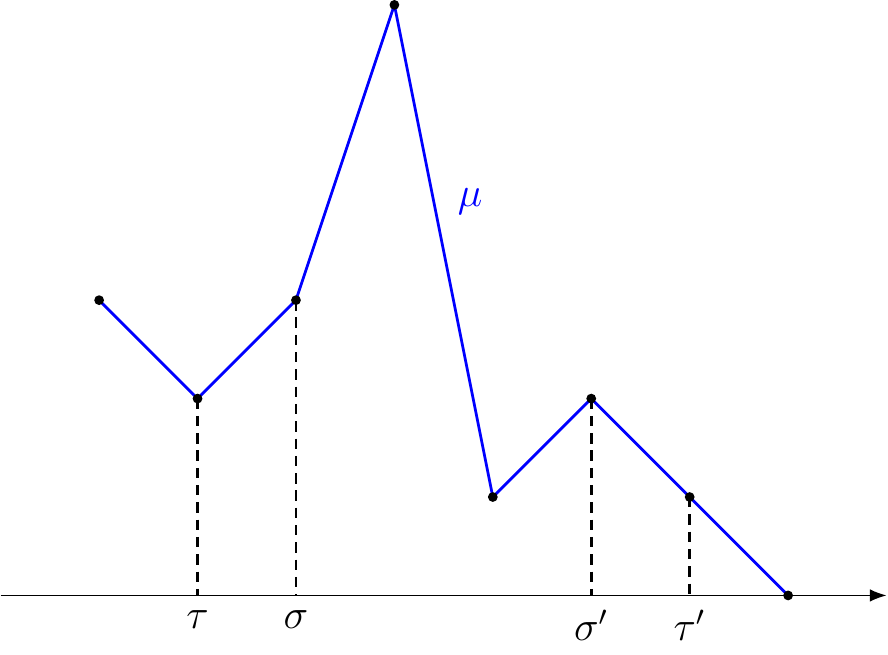}
  \end{subfigure}\hfill%
  \begin{subfigure}{0.48\linewidth}
    \includegraphics[width=\linewidth]{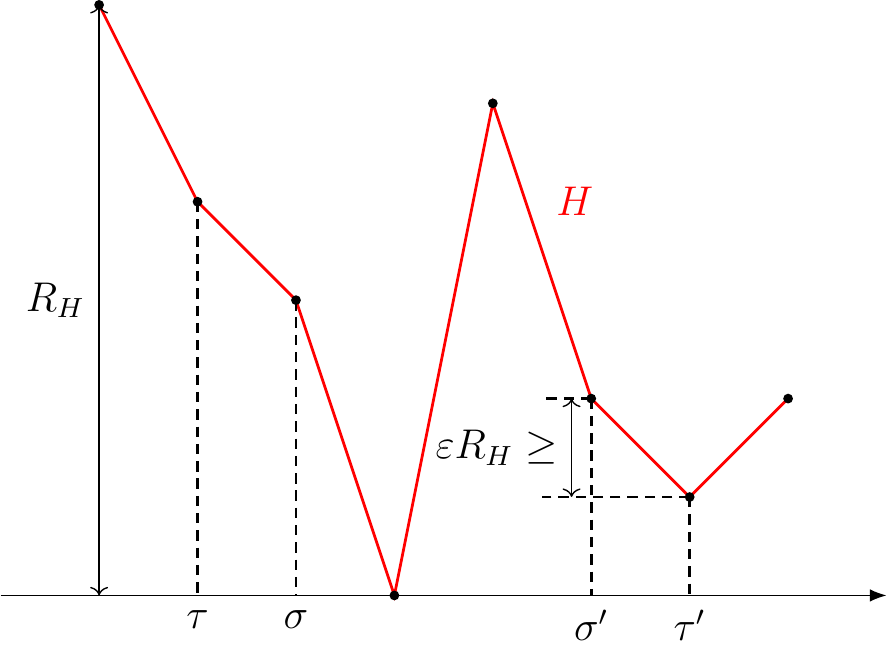}
  \end{subfigure}
  \caption{\label{fig:order-preservation}
Explanation of \eqref{eq:closedef} in two cases: either 
$\mu(\sigmavec)\ge\mu(\tauvec)$ and $H(\sigmavec)\le H(\tauvec)$ (the left 
pairs), or $\mu(\sigmavec')\ge\mu(\tauvec')$ and 
$H(\tauvec')<H(\sigmavec')\le H(\tauvec')+\vep R_H$ (the right pairs).
}
\end{figure}

Let $\vep=0.01$, for example.  
Then \refeq{closedef} means that the energy 
levels of two spin configurations respect their ordering in $\mu$ up to an 
error of 1\% of the range of $H$ (see Figure~\ref{fig:order-preservation}).  
We note that the latter inequality in \refeq{closedef} is equivalent to 
\begin{align}\lbeq{closerewr}
\piG_\beta(\sigmavec)\ge\piG_\beta(\tauvec)e^{-\beta\vep R_H},
\end{align}
for any $\beta$.

\begin{prp}\label{prp:close}
Let
\begin{align}\lbeq{vdef}
v=\sum_{\{x,y\}\in E}J_{x,y}^2+\sum_{x\in V}h_x^2,
\end{align}
and recall \refeq{Kdef} for the definition of $\bar K$.  Then, 
$\piSCA_{\beta,q}$ is $\vep$-close to the target Gibbs $\piG_\beta$ in 
the sense of order-preservation if
\begin{align}\lbeq{lowerbd}
2q\ge\log|V|+\beta\bar K-\log\frac{\vep\sqrt v}{2\bar K}.
\end{align}
\end{prp}

\Proof{Proof.}
First, by \refeq{sigmavec^I}, \refeq{cavity} and \refeq{symmetry}, we can 
rewrite $\wSCA_{\beta,q}(\sigmavec)/\wG_\beta(\sigmavec)$ as
\begin{align}
\frac{\wSCA_{\beta,q}(\sigmavec)}{\wG_\beta(\sigmavec)}&=\sum_{\tauvec}\exp
 \bigg(\beta\big(\tilde H(\sigmavec,\sigmavec)-\tilde H(\sigmavec,\tauvec)\big)
 +q\sum_{x\in V}\sigma_x\tau_x\bigg)\nn\\
&=e^{|V|q}\sum_{I\subset V}\sum_{\tauvec=\sigmavec^I}\exp\bigg(\frac\beta2
 \sum_{x\in V}\tilde h_x(\sigmavec)(\tau_x-\sigma_x)+q\sum_{x\in V}
 (\sigma_x\tau_x-1)\bigg)\nn\\
&=e^{|V|q}\sum_{I\subset V}\prod_{x\in I}e^{-\beta\tilde h_x(\sigmavec)\sigma_x
 -2q}\nn
=e^{|V|q}\prod_{x\in V}\big(1+\delta\phi_x(\sigmavec)\big),
\end{align}
where
\begin{align}
\delta=e^{-2q},&&
\phi_x(\sigmavec)=e^{-\beta\tilde h_x(\sigmavec)\sigma_x}.
\end{align}
As a result, 
\begin{align}
\frac{\piSCA_{\beta,q}(\sigmavec)}{\piSCA_{\beta,q}(\tauvec)}
 =\frac{\wSCA_{\beta,q}(\sigmavec)}{\wSCA_{\beta,q}(\tauvec)}
 =\frac{\wG_\beta(\sigmavec)}{\wG_\beta(\tauvec)}\prod_{x\in V}
 \frac{1+\delta\phi_x(\sigmavec)}{1+\delta\phi_x(\tauvec)}.
\end{align}

Now we suppose $\piSCA_{\beta,q}(\sigmavec)\ge\piSCA_{\beta,q}(\tauvec)$ 
(cf., the former inequality in \refeq{closedef}), so that
\begin{align}
\frac{\wG_\beta(\sigmavec)}{\wG_\beta(\tauvec)}\ge\prod_{x\in V}\frac{1+\delta
 \phi_x(\tauvec)}{1+\delta\phi_x(\sigmavec)}.
\end{align}
To prove the latter inequality in \refeq{closedef} (hence \refeq{closerewr}), 
it suffices to show 
\begin{align}
\prod_{x\in V}\frac{1+\delta\phi_x(\tauvec)}{1+\delta\phi_x(\sigmavec)}\ge
 e^{-\beta\vep R_H},
\end{align}
or equivalently 
\begin{align}\lbeq{vep-lb1}
\beta\vep R_H\ge\sum_{x\in V}\log\underbrace{\frac{1+\delta\phi_x(\sigmavec)}{1
 +\delta\phi_x(\tauvec)}}_{>0}=\sum_{x\in V}\log\bigg(1+\underbrace{\delta
 \frac{\phi_x(\sigmavec)-\phi_x(\tauvec)}{1+\delta\phi_x(\tauvec)}}_{>-1}\bigg).
\end{align}

Next we claim that 
\begin{align}\lbeq{vep-lb3}
R_H\ge\sqrt v,
\end{align}
where we recall \refeq{vdef} for the definition of $v$.  
This is due to the fact that 
\begin{align}
\mE_\mu\big[(H-\mE_\mu[H])^2\big]\le R_H^2
\end{align}
holds for any probability measure $\mu$ on $\{\pm1\}^V$, where $\mE_\mu$ 
denotes the expectation against $\mu$.  If $\mu$ is the uniform distribution 
over $\{\pm1\}^V$, then 
\begin{align}
\mE_\mu[H]&=\frac1{2^{|V|}}\sum_{\sigmavec\in\{\pm1\}^V}H(\sigmavec)=0,\\
\mE_\mu[H^2]&=\frac1{2^{|V|}}\sum_{\sigmavec\in\{\pm1\}^V}H(\sigmavec)^2=v,
\end{align}
which yields \refeq{vep-lb3}.

Next we claim that 
\begin{align}\lbeq{vep-lb2}
\sum_{x\in V}\log\bigg(1+\delta\frac{\phi_x(\sigmavec)-\phi_x(\tauvec)}{1+\delta
 \phi_x(\tauvec)}\bigg)\le2|V|\beta\delta\bar K\,e^{\beta\bar K},
\end{align}
where we recall \refeq{Kdef} for the definition of $\bar K$.  
To show this, the first step is to use $\log(1+a)\le a\le|a|$ for all $a>-1$.  
Then we obtain
\begin{align}
\log\bigg(1+\delta\frac{\phi_x(\sigmavec)-\phi_x(\tauvec)}{1+\delta\phi_x
 (\tauvec)}\bigg)\le\bigg|\delta\frac{\phi_x(\sigmavec)-\phi_x(\tauvec)}{1
 +\delta\phi_x(\tauvec)}\bigg|\le\delta|\phi_x(\sigmavec)-\phi_x(\tauvec)|.
\end{align}
Since $|e^X-e^Y|\le|X-Y|(e^X\vee e^Y)$ for any real $X,Y$, it is further 
bounded above by
\begin{align}\lbeq{logphi-lb1}
\beta\delta\,|\tilde h_x(\sigmavec)\sigma_x-\tilde h_x(\tauvec)\tau_x|\big(
 \phi_x(\sigmavec)\vee\phi_x(\tauvec)\big)\le2\beta\delta\bar K\,e^{\beta
 \bar K},
\end{align}
which yields \refeq{vep-lb2}.

Combining the above two claims \refeq{vep-lb3} and \refeq{vep-lb2}, we arrive 
at a sufficient condition for \refeq{vep-lb1}:
\begin{align}
\vep\sqrt v\ge2|V|\delta\bar K\,e^{\beta\bar K},
\end{align}
which is equivalent to \refeq{lowerbd}.
\QED

\section{Concluding remark}\label{s:concluding}
Propositions~\ref{prp:morespins} and \ref{prp:close} imply the following.

\begin{thm}\label{thm:summary}
If
\begin{align}\lbeq{betabd}
\beta\le\frac1{2\bar K}\log\frac{\vep\sqrt v}{2\bar K},
\end{align}
then there is a pinning parameter $q$ for which the following both hold at the 
same time:
\begin{itemize}
\item
SCA is faster than Glauber in terms of the expected number of spin-flips per 
update, i.e., 
$\EG_\beta[|D_{\sigmavec,X_{\sigmavec}}|]\le\ESCA_{\beta,q}[|D_{\sigmavec,
X_{\sigmavec}}|]$ for any $\sigmavec\in\{\pm1\}^V$.
\item
The SCA equilibrium $\piSCA_{\beta,q}$ is $\vep$-close to the target Gibbs 
$\piG_\beta$ in the sense of order-preservation.  In particular, if 
$\piSCA_{\beta,q}$ takes on its highest point at $\sigmavec$, then 
$H(\sigmavec)\le H(\tauvec)+\vep R_H$ holds for all $\tauvec\in\{\pm1\}^V$.
\end{itemize}
\end{thm}

The main message of the above theorem is that, to ensure the above two 
properties simultaneously, we have to keep the temperature sufficiently high.  
However, if the temperature is too high (e.g., $\beta=0$), then the equilibrium 
distributions are almost identical to the uniform distribution over 
$\{\pm1\}^V$, which is not at all efficient to search for the ground states.  
In this respect, we want to take the largest possible value of $\beta$.  

The upper bound in \refeq{betabd} depends on the coupling coefficients 
$\{J_{x,y}\}_{\{x,y\}\in E}$ and $\{h_x\}_{x\in V}$ as well as the graph 
structure $G=(V,E)$.  Suppose that the coupling coefficients 
are uniformly bounded (e.g., $J_{x,y},h_x\in\{\pm1\}$).  If $G$ is a regular 
graph with bounded degree (e.g., a subset of the square lattice), then the 
bound in \refeq{betabd} is of order $\log|V|$, which is better than 
the bound obtained by Dobrushin's condition (e.g., \cite{dss12}) in a large 
volume, and we may find a good candidate for the ground states of a large 
system via a sufficiently slow cooling method $\beta\to\infty$, such as 
simulated annealing \cite{c85,kgv83}.  On the other hand, if $G$ is a complete 
graph, then the bound in \refeq{betabd} is finite (or, even worse, may be 
negative) and therefore the two properties in Theorem~\ref{thm:summary} may 
not hold at the same time under any cooling methods (see \cite{st18} for 
application of SCA to the SK spin glass).

Instead of taking $\beta$ to infinity, we should keep $\beta$ bounded as 
in \refeq{betabd}.  Since $\piSCA_{\beta,q}$ (supposed to be close to 
$\piG_\beta$) is not concentrated on the ground states, the chance of a single 
MCMC experiment hitting a ground states is slim.  Instead, we may have to run 
$N$ independent MCMCs $\{X_j\}_{j=1}^N$ whose common law is 
$U*(\PSCA_{\beta,q})^{*n}$, where $U$ is the uniform distribution over 
$\{\pm1\}^V$, and construct a profile $\frac1N\sum_{j=1}^N\delta_{X_j}$ 
that is supposed to be close to $\piSCA_{\beta,q}$ (hence to $\piG_\beta$) if 
$n$ and $N$ are sufficiently large.  Then, a ground state may be captured as 
$\arg\max\frac1N\sum_{j=1}^N\delta_{X_j}$.  How large $n$ has to be depends 
on the mixing time, which is under investigation with Bruno Kimura.  On the 
other hand, how large $N$ has to be depends on the convergence rate in the 
law of large numbers (LLN).  Sanov's theorem \cite{s57} gives the exact 
exponential rate of convergence, but the multiplicative term could be huge (in 
powers of $N$) depending on the degree of degeneracy.  This is partly due to 
the fact that the number $2^{|V|}$ of spin configurations is way larger than 
the number $N$ of experiments, which is not in the regime of LLN.  We are 
currently seeking to apply a theory of big data: high-dimensional statistics 
\cite{g15}.

\section*{Acknowledgements}
This work was supported by JST CREST Grant Number JP22180021, Japan.  We would 
like to thank Takashi Takemoto of Hitachi, Ltd., for providing us with a 
stimulating platform for the weekly meeting at Global Research Center for Food 
\& Medical Innovation (FMI) of Hokkaido University.  We would also like to 
thank Mertig Normann of Hitachi, Ltd., and Hiroshi Teramoto of Research 
Institute for Electronic Science (RIES), as well as Ko Fujisawa, Masamitsu 
Aoki, Yuki Ueda, Hisayoshi Toyokawa, Naomichi Nakajima and Shinpei Makita of 
Mathematics Department, for valuable comments and encouragement at the 
aforementioned meetings at FMI.  We also extend our thanks to Yuki Chino, 
Eric Endo and Bruno Kimura for intensive discussion during their visits to 
Hokkaido University from February 12 through 23, 2019.

\end{document}